\documentclass[12pt]{article}
\usepackage{latexsym, amscd, %graphicx, 
amssymb}
    \title{{\bf  Riemann surfaces with boundaries and the theory of 
vertex operator algebras}}
    \author{Yi-Zhi Huang}
    \date{}
    
\begin{document}
    \bibliographystyle{alpha}
    \maketitle

   \newtheorem{thm}{Theorem}[section]
\newtheorem{defn}[thm]{Definition}
\newtheorem{prop}[thm]{Proposition}
\newtheorem{cor}[thm]{Corollary}
\newtheorem{rema}[thm]{Remark}
\newtheorem{lemma}[thm]{Lemma}
\newtheorem{app}[thm]{Application}
\newtheorem{prob}[thm]{Problem}
\newtheorem{conv}[thm]{Convention}
\newtheorem{conj}[thm]{Conjecture}
\newcommand{\halmos}{\rule{1ex}{1.4ex}}
\newcommand{\pfbox}{\hspace*{\fill}\mbox{$\halmos$}}

 \newcommand{\nn}{\nonumber \\}
 
	\newcommand{\nno}{\nonumber}
	\newcommand{\lbar}{\bigg\vert}
\newcommand{\mbar}{\mbox{\large $\vert$}}
	\newcommand{\p}{\partial}
	\newcommand{\dps}{\displaystyle}
	\newcommand{\bra}{\langle}
	\newcommand{\ket}{\rangle}
 \newcommand{\res}{\mbox{\rm Res}}
\renewcommand{\hom}{\mbox{\rm Hom}}
\newcommand{\hol}{\mbox{\rm Hol}}
\newcommand{\dt}{\mbox{\rm Det}}
\newcommand{\edo}{\mbox{\rm End}\;}
 \newcommand{\pf}{{\it Proof.}\hspace{2ex}}
 \newcommand{\epf}{\hspace*{\fill}\mbox{$\halmos$}}
 \newcommand{\epfv}{\hspace*{\fill}\mbox{$\halmos$}\vspace{1em}}
 \newcommand{\epfe}{\hspace{2em}\halmos}
\newcommand{\nord}{\mbox{\scriptsize ${\circ\atop\circ}$}}
\newcommand{\wt}{\mbox{\rm wt}\ }
\newcommand{\tr}{\mbox{\rm Tr}}
\newcommand{\swt}{\mbox{\rm {\scriptsize wt}}\ }
\newcommand{\clr}{\mbox{\rm clr}\ }

\begin{abstract} 
The connection between Riemann surfaces with boundaries
and the theory of vertex operator algebras is discussed in the framework
of conformal field theories defined by Kontsevich and Segal and in
the framework of their generalizations in open string theory and
boundary conformal field theory.  We present some results, problems,
conjectures, their conceptual implications and meanings in a program to
construct these theories from representations of vertex operator
algebras. 
\end{abstract}

\renewcommand{\theequation}{\thesection.\arabic{equation}}
\renewcommand{\thethm}{\thesection.\arabic{thm}}
\setcounter{equation}{0}
\setcounter{thm}{0}
\setcounter{section}{0}
\section{Introduction}

Quantum field theory is one of the greatest gifts that physicists have
brought to mathematicians. In physics, the quantum-field-theoretic
models of electro-magnetism, weak and strong interactions are among the
most successful theories in physics. Besides its direct applications in
building realistic physical models, quantum field theory is also an
important ingredient and a powerful tool in string theory or M theory,
which is being developed by physicists as a promising candidate of a unified
theory of all interactions including gravity. To mathematicians, the
most surprising thing is that ideas and intuition in quantum field
theory have been used successfully by physicists to make deep
mathematical conjectures and to discover the unity and beauty of different
branches of mathematics in connection with the real world.   The quantum 
invariants of knots and three-dimensional manifolds,
Verlinde formula, mirror symmetry and Seiberg-Witten theory are among
the most famous examples. The results predicted by these physical ideas
and intuition also suggest that many seemingly-unrelated mathematical
branches are in fact different aspects of a certain
yet-to-be-constructed unified theory. The success of physical ideas and
intuition in mathematics has shown that, no matter how abstract it is, our
mathematical theory is deeply related to the physical world in which we
live. 

One of the most challenging problem for mathematicians is to understand
precisely the deep ideas and conjectures obtained by physicists using
quantum field theories. In particular, quantum field theories must be
formulated precisely and constructed mathematically. The modern approach
to quantum field theory used by physicists is the one based on path
integrals. Though path integrals are in general not defined
mathematically, their conjectured properties can be extracted to give
axiomatic definitions of various notions of quantum field theory. When a
definition of quantum field theory is given mathematically, a
well-formulated problem is how to construct such a theory. The most
successful such theories are topological quantum field theories, since the
state spaces of these theories are typically finite-dimensional and
consequently  many
existing mathematical theories can readily be applied to the
construction. Outside the realm of topological  theories,
quantum field theories seem to be very difficult to construct and study
because any nontrivial nontopological theory  must have an
infinite-dimensional state space. 

Fortunately there is one simplest class of nontopological quantum field
theories---two-dimensional conformal field theories---which has been
studied extensively by physicists and mathematicians using various
approaches.  
It is expected that a complete mathematical understanding
of these conformal field theories will not only solve problems related
to these theories, but will also give hints and insights to other classes
of quantum field theories.  

These theories arose in both condensed matter physics and string theory.
The systematic development of (two-dimensional) conformal field theories
in physics started with the seminal paper by Belavin, Polyakov and
Zamolodchikov \cite{BPZ} and the study of their geometry was initiated
by Friedan and Shenker in \cite{FS}.  Mathematically, around the same
time, using the theory of vertex operators developed in the early days
of string theory and in the representation theory of affine Lie
algebras, Frenkel, Lepowsky and Meurman \cite{FLM1} gave a construction
of what they called the ``moonshine module,'' an infinite-dimensional
representation of the Fischer-Griess Monster finite simple group and,
based on his insight in representations of affine Lie algebras and the
moonshine module, Borcherds \cite{B} introduced the notion of vertex
algebra. The theory of vertex operator algebras was further developed by
Frenkel, Lepowsky and Meurman in \cite{FLM} and, in the same work, the
statement by Borcherds in \cite{B} that the moonshine module has a
structure of a vertex operator algebra was proved.  Motivated by the
path integral formulation of string theory, I. Frenkel initiated a
program to construct geometric conformal field theories in a suitable
sense in 1986. Also around the same time, Y.~Manin \cite{Man}
pointed out that the moduli space of algebraic curves of all genera
should play the same role in the representation theory of the Virasoro
algebra as the space $G/P$ ($P$ is a parabolic subgroup) in the
representation theory of a semisimple Lie group $G$. Soon M.~Kontsevich
\cite{K} and A.~Beilinson and V.~Schechtman \cite{BS} found the
relationship between the Virasoro algebra and the determinant line
bundles over the moduli spaces of curves with punctures.

In 1987, Kontsevich and Segal
independently gave a precise definition of conformal field theory using
the properties of path integrals as axioms (see \cite{S1}). 
To explain the rich structure of chiral parts of conformal field theories, 
Segal in \cite{S2} and \cite{S3} further introduced the notions of
modular functor and weakly conformal field theory and sketched how to
obtain a conformal field theory from a suitable weakly conformal field
theory. In \cite{V},  C.~Vafa also gave, on a
physical level of rigor, a formulation of conformal field theories using
Riemann surfaces with punctures and local coordinates vanishing at these
punctures.

Starting from Segal's axioms and some additional properties of
rational conformal field theories, Moore and Seiberg \cite{MS1} \cite{MS2}
constructed modular tensor categories and proved the
Verlinde conjecture, which states that the fusion rules are diagonalized
simultaneously by the modular transformation corresponding to
$\tau\mapsto -1/\tau$.  Around the same time, 
Tsuchiya, Ueno and Yamada \cite{TUY}  constructed the 
the algebro-geometric parts of conformal field theories corresponding to 
the Wess-Zumino-Novikov-Witten models. The 
algebro-geometric parts of conformal field theories corresponding to 
the minimal models were later constructed by Beilinson, Feigin and Mazur
\cite{BFM}. Though these algebro-geometric parts 
of conformal field theories do not give (weakly) 
conformal field theories
in the sense of Segal, it seems that they
are needed, explicitly or implicitly, in any future construction of 
full theories.

The definition of conformal field theory by Kontsevich and Segal is
based on Riemann surfaces with boundaries. It was first observed by
I. Frenkel that vertex operators, when modified slightly,
actually correspond to
the unit disk with two smaller disks removed. Based on a thorough study
of the sewing operation for the infinite-dimensional moduli space of
spheres with punctures and local coordinates and the determinant 
line bundle over this moduli space, the author in \cite{H1},
\cite{H2} and \cite{H3}
gave a geometric definition of vertex operator algebra and 
proved that the geometric and the algebraic definitions
are equivalent. In \cite{Z}, Y. Zhu proved that for a suitable vertex 
operator algebra, the $q$-traces of products of vertex operators associated 
to modules are modular invariant in a certain sense. 
The moduli space mentioned above 
contains the moduli space of genus-zero Riemann
surfaces with boundaries and analytic boundary parametrizations, 
and therefore we
expect that vertex operator algebras will play a crucial role in the
construction of conformal field theories.

Recently, in addition to the continuing development of conformal field
theories, several new directions have attracted attentions. For example,
boundary conformal field theories first developed by Cardy in \cite{C1}
and \cite{C2} play a fundamental role in many problems in condensed
matter physics, and they have also become one of the main tools in the
study of $D$-branes, which are nonperturbative objects in string theory.
In the framework of topological field theories, boundary topological
field theories (open-closed topological field theories) have been
studied in detail by Lazaroiu \cite{L} and by Moore and Segal \cite{Mo}
\cite{S4}. Some part of the analogue in the conformal case of the work
in \cite{L}, \cite{Mo} and \cite{S4} has been done mathematically by
Felder, Fr\"{o}hlich, Fuchs and Schweigert \cite{FFFS} and by Fuchs,
Runkel and Schweigert \cite{FRS1} (see also \cite{FRS2}).  But in this
conformal case, even boundary conformal field theories (``open-closed
conformal field theories'') on disks with parametrized boundary segments
still need to be fully constructed and studied and the mathematics involved is
obviously very hard but also very deep.  On the other hand, many of the
ideas, mathematical tools and structures used in the study of conformal
field theories can be adapted to the study of open-closed conformal
field theories and $D$-branes, which can be viewed as substructures in
open-closed conformal field theories. Besides the obvious problem of
constructing and classifying open-closed conformal field theories, the
study of $D$-branes and their possible use in geometry have led to
exciting and interesting mathematical problems. 

Logarithmic conformal field theories initiated by Gurarie in \cite{G}
provide another new direction.  These theories, which describe
disordered systems in condensed matter physics, also occur naturally in
the mathematical study of conformal field theories. In \cite{M}, Milas
formulated and studied some of the basic ingredients of logarithmic
conformal field theories in terms of representations of vertex operator
algebras.  One expects that there should also be a geometric formulation
of logarithmic conformal field theory. The main mathematical problems
would consist of the construction and classification of such theories.
The representation theory of the Virasoro algebra with central charge 0
provides a crucial example of such a theory. We shall not discuss
logarithmic conformal field theories in this paper. Instead, we refer
the reader to the papers mentioned above and many other papers by
physicists (for example, the expositions \cite{Ga} by Gaberdiel, \cite{F}
by Flohr and 
the references there).

In the present paper, we shall discuss some results, problems,
conjectures, their conceptual implications and meanings in a program 
to construct conformal field theories and their generalizations 
from representations of vertex operator algebras. In the next section,
we recall a definition of (closed) 
conformal field theory by Kontsevich and Segal
and other important notions introduced by Segal. We discuss the program
of constructing conformal field theories in Section 3. In Section 4,
we discuss open-closed conformal field theories which incorporate the
axiomatic properties of conformal field theories and boundary
conformal field theories.

\paragraph{Acknowledgment} This research is supported in part by NSF
grant DMS-0070800. The present paper was an expanded and updated version
of the notes of the author's talk at the workshop, written when the
author was visiting the Department of Mathematics at University of
Virginia. The author would like to thank David Radnell for carefully
reading a draft of the paper and J\"{u}rgen Fuchs, Liang Kong, Antun Milas, 
David Radnell, Christoph Schweigert and Jim Stasheff for comments 
and/or discussions.

\renewcommand{\theequation}{\thesection.\arabic{equation}}
\renewcommand{\thethm}{\thesection.\arabic{thm}}
\setcounter{equation}{0}
\setcounter{thm}{0}

\section{Conformal field theories}

In this section, we recall the definition of conformal field theory
in the sense of Kontsevich and Segal and other notions introduced by Segal.

First, we recall the definition of conformal field theory.  For more
details, see \cite{S1}, \cite{S2} and \cite{S3}.  Consider the following
symmetric monoidal category $\mathcal{B}$ constructed geometrically:
Objects of $\mathcal{B}$ are finite ordered sets of copies of the unit
circle $S^{1}$. Given two objects, morphisms from one object to another
are conformal equivalence classes of Riemann surfaces (including
degenerate ones, e.g., circles, and possibly disconnected) with oriented
and ordered boundary components together with analytic parametrizations
of these components such that the copies of $S^{1}$ in the domain and
codomain parametrize the negatively oriented and positively oriented
boundary components, respectively.  For an object containing $n$ unit
circles, the identity on it is the degenerate surface given by the $n$
unit circles with the trivial parametrizations of the boundary
components.  Given two composable morphisms, we can compose them using
the boundary parametrizations in the obvious way.  It is easy to see
that the composition satisfies the associativity and any morphism
composed with an identity is equal to itself. Thus we have a category.
This category has a symmetric monoidal category structure defined by
disjoint unions of objects and morphisms. We shall use $[\Sigma]$ to denote
the conformal equivalence class of a Riemann surfaces  $\Sigma$ 
with oriented,
ordered and parametrized boundary components.

We also have a symmetric tensor category $\mathcal{T}$ of complete
locally convex topological vector spaces over $\mathbb{C}$
with nondegenerate bilinear
forms.  A projective functor from $\mathcal{B}$ to $\mathcal{T}$ is a
functor from $\mathcal{B}$ to the projective category of $\mathcal{T}$
which has the same objects as $\mathcal{T}$ but the morphisms are 
one-dimensional spaces of morphisms of $\mathcal{T}$.  Note that for any
functor or projective functor $\Phi$ of monoidal categories from
$\mathcal{B}$ to $\mathcal{T}$, if $\Phi(S^{1})=H$, then the image under
$\Phi$ of the object containing $n$ ordered copies of $S^{1}$ must be
$H^{\otimes n}$.

A {\it conformal field theory} (or {\it closed conformal field theory})
is a projective functor $\Phi$ from the
symmetric monoidal category $\mathcal{B}$ to the symmetric tensor
category $\mathcal{T}$ satisfying the following axioms: (i) Let $[\Sigma]$
be a morphism in $\mathcal{B}$ from $m$ ordered copies of $S^{1}$ to $n$
ordered copies of $S^{1}$. Let $[\Sigma_{\widehat{\widehat{i, j}}}]$ 
be the morphism
from $m-1$ copies of $S^{1}$ to $n-1$ copies of $S^{1}$ obtained from
$[\Sigma]$ by identifying the boundary component of $\Sigma$
parametrized by the $i$-th copy of $S^{1}$ in the domain of
$[\Sigma]$ with the boundary component of $\Sigma$
parametrized by the $j$-th copy of $S^{1}$ in the codomain of $[\Sigma]$.
(Note that the two copies of $S^{1}$ identified might or might not be 
on a same connected component of $\Sigma$.)
Then the trace between the $i$-th tensor factor of the domain and the
$j$-th tensor factor of the codomain of $\Phi([\Sigma])$ exists and is
equal to $\Phi([\Sigma_{\widehat{\widehat{i,j}}}])$.  (ii) Let $[\Sigma]$ be a
morphism in $\mathcal{B}$ from $m$ ordered copies of $S^{1}$ to $n$
ordered copies of $S^{1}$. Let $[\Sigma_{i\to n+1}]$ be the morphism from
the set of $m-1$ ordered copies of $S^{1}$ to the set of $n+1$ ordered
copies of $S^{1}$ obtained by changing the $i$-th copy of $S^{1}$ of the
domain of $[\Sigma]$ to the $n+1$-st copy of $S^{1}$ of the codomain of
$[\Sigma_{i\to n+1}]$. Then
$\Phi([\Sigma])$ and $\Phi([\Sigma_{i\to n+1}])$ are related by the map from
$\hom(H^{\otimes m}, H^{\otimes n})$ to $\hom(H^{\otimes m-1},
H^{\otimes (n+1)})$ obtained using the map $H\to H^{*}$ corresponding to
the bilinear form $(\cdot, \cdot)$.

A {\it real conformal field theory}  is a conformal field theory 
together with an anti-linear involution $\theta$ from $H$ to itself
satisfying the following additional axiom: Let $[\Sigma]$ be a morphism in 
$\mathcal{B}$ from $m$ ordered copies of
$S^{1}$ to $n$ ordered copies of $S^{1}$ and $[\overline{\Sigma}]$ the 
morphism in $\mathcal{B}$ from $n$ ordered copies of
$S^{1}$ to $m$ ordered copies of $S^{1}$ obtained by taking the complex 
conjugate complex structure of the one on $[\Sigma]$ 
(note that the orientations
of the boundary components are reversed). Then 
$\Phi([\overline{\Sigma}])=\theta^{\otimes m}\circ \Phi^{*}([\Sigma])\circ
(\theta^{-1})^{\otimes n}$ where $\Phi^{*}([\Sigma])$ is the adjoint 
of $\Phi([\Sigma])$.

The definitions above do not reveal many important ingredients in
concrete models. In particular, they do not give the detailed structure
of chiral and anti-chiral parts of conformal field theories, that is,
parts of conformal field theories depending on the moduli space
parameters analytically and anti-analytically.  It is known that
meromorphic fields in a conformal field theory form a vertex operator
algebra. The representations of this vertex operator algebra 
form the chiral parts of the theory.
Therefore to construct conformal field theories from vertex
operator algebras, it is necessary to study first chiral and anti-chiral
parts of conformal field theories.  Axiomatically, chiral and
anti-chiral parts of conformal field theories are weakly conformal field
theories defined by G. Segal in \cite{S2} and \cite{S3} and are
generalizations of conformal field theories. 

To describe weakly conformal field theories, we first need to describe
modular functors. Instead of Riemann surfaces with parametrized boundary
components, we need {\it Riemann surfaces with labeled, oriented and
parametrized boundaries}, which are Riemann surfaces with oriented and
parametrized boundaries and an assignment of an element of a fixed set
$\mathcal{A}$ to each boundary component. The set $\mathcal{A}$ is
typically the set of equivalence classes of irreducible modules for a
vertex operator algebra.  We consider a category whose objects are
conformal equivalence classes of Riemann surfaces with labeled, oriented
and parametrized boundaries and whose morphisms are given by the sewing
operation, that is, if one such equivalent class can be obtained from
another using the sewing operation, then the procedure of obtaining the
second surface from the first one is a morphism. We use $[\Sigma]$ to denote 
the conformal equivalence class of the surface $\Sigma$. 

A {\it modular functor} is a functor $E$ from the above category of
Riemann surfaces with labeled, oriented and parametrized boundaries to the
category of finite-dimensional vector spaces over $\mathbb{C}$ 
satisfying the following
conditions: (i) $E([\Sigma_{1} \sqcup \Sigma_{2}])$ is naturally
isomorphic to $E([\Sigma_{1}])\otimes E([\Sigma_{2}])$.  (ii) If $\Sigma$ is
obtained from another surface $\Sigma_{a}$ by sewing two boundary
components with opposite orientations but with the same label $a\in
\mathcal{A}$ of $\Sigma_{a}$, then $E([\Sigma])$ is naturally isomorphic
to $\oplus_{b\in \mathcal{A}}E([\Sigma_{b}])$ where for $b\ne a$,
$\Sigma_{b}$ is the surface obtained from $\Sigma_{a}$ by changing the
label $a$ to $b$ on the boundary components to be sewn. (iii) $\dim
E([S^{2}])=1$. (iv) $E([\Sigma])$ depends on $\Sigma$ holomorphically.

From a modular functor $E$, we can construct a symmetric monoidal
category $\mathcal{B}_{E}$ extending the category $\mathcal{B}$ as
follows: Objects of $\mathcal{B}_{E}$ are ordered sets of pairs of the
form of a copy of the unit circle $S^{1}$ and an element of the set
$\mathcal{A}$. Morphisms of $\mathcal{B}_{E}$ are pairs of the form of an
equivalence class $[\Sigma]$ of 
Riemann surface  with labeled, oriented and parametrized
boundaries and the vector space $E([\Sigma])$, such that the labels of the
boundary components of $\Sigma$ match with the labels of the copies of
$S^{1}$ in the domain and codomain.  Let $E$ be a modular functor
labeled by $\mathcal{A}$.  A {\it weakly conformal field theory over
$E$} is a functor $\Phi$ from the symmetric monoidal category
$\mathcal{B}_{E}$ to the symmetric tensor category $\mathcal{T}$ 
satisfying the obvious axioms similar to those for conformal 
field theories.

\renewcommand{\theequation}{\thesection.\arabic{equation}}
\renewcommand{\thethm}{\thesection.\arabic{thm}}
\setcounter{equation}{0}
\setcounter{thm}{0}

\section{Vertex operator algebras and conformal field theories}

The definitions of the notions in the preceding section of conformal
field theory, modular functor and weakly conformal field theory are
simple, beautiful and conceptually satisfactory. The obvious first
problem is the following:

\begin{prob}
Construct conformal field theories, real conformal field theories,
modular functors and weakly conformal field theories.
\end{prob}

Since the definition of (weakly) conformal field theory involves
algebra, geometry and analysis, it is not surprising that the
construction of examples is very hard. Up to now, for only the free
fermion theories has a complete construction been sketched, by Segal in
an unpublished manuscript; this construction has been further clarified
recently by I.~Kriz \cite{Kr}.  Also the method used for free fermion
theories does not generalize to other cases.  In \cite{S2}, Segal
described a construction of conformal field theories from weakly
conformal field theories satisfying a "unitarity" condition. Thus the
main hard problem is actually the construction of modular functors and
weakly conformal field theories.  In this section, we discuss results
and problems in a program to construct weakly conformal field theories
in the sense of Segal from representations of suitable vertex operator
algebras.

Since Riemann surfaces with labeled, oriented
and parametrized boundaries can be
decomposed into genus-zero Riemann surfaces with labeled, oriented and
parametrized boundaries, the first step is to construct the genus-zero
weakly conformal field theories, that is, to construct spaces associated
to the unit circle with labels, to construct maps associated to genus-zero
surfaces and to prove the axioms which make sense for genus-zero
surfaces. This step has been worked out by the author using the
theory of intertwining operator algebras developed by the author and 
the tensor product theory for modules for a vertex operator algebra of
Lepowsky and the author.  See \cite{H1}--\cite{H11} and
\cite{HL3}--\cite{HL6} for the construction and the theories used.

For basic notions in the theory of vertex operator algebras, we 
refer the reader to \cite{FLM} and \cite{FHL}. 
We need the following notions to state the main result: Let $V$ be a vertex 
operator algebra and $W$ a $V$-module. Let $C_{1}(W)$
be the subspace of $W$ spanned by elements of the form $u_{-1}w$ for
$u\in V_{+}=\coprod_{n>0}V_{(n)}$ and $w\in W$. If $\dim
W/C_{1}(W)<\infty$, we say that $W$ is {\it $C_{1}$-cofinite} or $W$
satisfies the {\it $C_{1}$-cofiniteness condition}. 
A {\it generalized $V$-module} is a pair 
$(W, Y)$ satisfying all the 
conditions for a $V$-module except for
the two grading-restriction conditions.

The main result in the genus-zero case is the following:

\begin{thm}\label{genus-zero}
Let $V$ be a vertex operator algebra satisfying the 
following conditions:

\begin{enumerate}

\item  Every  generalized $V$-module is a direct sum of 
irreducible $V$-modules.

\item There are only finitely many inequivalent irreducible $V$-modules
and they are all $\mathbb{R}$-graded.  

\item Every
irreducible $V$-module satisfies the $C_{1}$-cofiniteness condition.

\end{enumerate}

Then  there is a natural structure of genus-zero weakly conformal
field theory with the set $\mathcal{A}$ of equivalence classes 
of irreducible $V$-modules as the set of labels and 
suitable locally convex topological completions $H^{a}$
of representatives of equivalence classes $a\in \mathcal{A}$
as the spaces corresponding to the unit circle with label $a$. 
\end{thm}

For vertex operator algebras associated to minimal models, WNZW
models, $N=1$  superconformal minimal models and $N=2$ 
superconformal unitary models, the conditions of Theorem 
\ref{genus-zero} are satisfied.

The second logical step is to construct genus-one  theories, 
that is, to construct maps associated to genus-one surfaces
and prove the axioms which make sense for genus-one surfaces. 
The first result in this step was obtained by Zhu \cite{Z}. To state this 
result, we need some further notions: 
Let $V$ be a vertex operator algebra and $C_{2}(V)$ is the subspace of 
$V$ spanned by elements of the form $u_{-2}v$ for
$u, v\in V$. Then we say that $V$ is {\it $C_{2}$-cofinite}
or satisfies the {\it $C_{2}$-cofiniteness
condition} if $V/C_{2}(V)$ is finite-dimensional. 
An {\it $\mathbb{N}$-gradable weak $V$-module} is a pair 
$(W=\coprod_{n\in \mathbb{N}}W_{(n)}, Y)$
satisfying all the conditions for a $V$-module except that the 
grading-restriction conditions are not required and the 
$L(0)$-grading condition is replaced by the condition that 
$u_{n}W_{(k)} \subset W_{(m-n-1+k)}$ for $u\in V_{(m)}$,
$n\in \mathbb{Z}$ and $k\in \mathbb{N}$. The following is 
Zhu's result (see \cite{Z}):

\begin{thm}\label{zhu}
Let $V$ be a vertex operator algebra of central charge $c$
satisfying the following 
conditions:

\begin{enumerate}

\item $V$ has no nonzero elements of negative weights.

\item Every $\mathbb{N}$-gradable $V$-module
is completely reducible.

\item There are only finitely
many irreducible $V$-modules.

\item  $V$ satisfies the 
$C_{2}$-cofiniteness condition.

\item $V$ as a module for the Virasoro algebra is generated
by lowest weight vectors.

\end{enumerate}

Let $\{M_{1}, \dots, M_{m}\}$ be a complete set of 
representatives of equivalence 
classes of (inequivalent) irreducible $V$-modules. Then for any
$n\in \mathbb{Z}_{+}$, we have:

\begin{enumerate}

\item For $i=1, \dots, m$, $u_{1}, \dots, u_{n}\in V$,
$$\tr_{M_{i}}Y(z_{1}^{L(0)}u_{1}, z_{1})
\cdots Y(z_{n}^{L(0)}u_{n}, z_{n})q^{L(0)-\frac{c}{24}}$$
is absolutely convergent when $1>|z_{1}|>\cdots>|z_{n}|>|q|>0$
and can be analytically extended to a meromorphic function
in the region
$1>|z_{1}|, \dots, |z_{n}|>|q|>0$, $z_{i}\ne z_{j}$ when $i\ne j$.

\item Let
$$S_{M_{i}}((u_{1}, z_{1}), \dots, (u_{n}, z_{n}), \tau)$$
for $i=1, \dots, m$ be the analytically extensions above
with $z_{j}$ replaced by $e^{2\pi iz_{j}}$ for $j=1, \dots, n$, and
$q$ replaced by $e^{2\pi \tau}$. 
For $i=1, \dots, m$, $u_{1}, \dots, u_{n}\in V$ and 
$$\left(\begin{array}{cc}
\alpha&\beta\\
\gamma&\delta
\end{array}\right)\in {\rm SL}(2, \mathbb{Z}),$$
$$S_{M_{i}}\left(\left(\left(\frac{1}{\gamma \tau+\delta}\right)^{L(0)}
u_{1}, \frac{z_{1}}{\gamma \tau+\delta}\right), 
\dots, \left(\left(\frac{1}{\gamma \tau+\delta}\right)^{L(0)}
u_{n}, \frac{z_{n}}{\gamma \tau+\delta}\right), \tau\right)$$
is a linear combination of 
$$S_{M_{j}}((u_{1}, z_{1}), \dots, (u_{n}, z_{n}), \tau)$$
for $j=1, \dots, m$. 

\end{enumerate}
\end{thm}

In \cite{H11},  it was proved that the conclusion of Theorem \ref{genus-zero}
is also true if for $n<0$, $V_{(n)}=0$  and $V_{(0)}=\mathbb{C}\mathbf{1}$
and Conditions 2-4 in Theorem \ref{zhu} hold.

By modifying Zhu's method, Dong, Li, Mason \cite{DLM} and Miyamoto
\cite{Mi1} \cite{Mi2} generalized Zhu's result above in a number of
directions. In particular, it was shown in \cite{DLM} that 
Condition 5 is not needed in Theorem \ref{zhu}.
But to construct genus-one theories completely, we need to
generalize Zhu's result to intertwining operator algebras. In this
general case, Zhu's method cannot be modified to study the $q$-traces
of products of more than one intertwining operators.  Recently, the
author has solved this problem.  Here we describe the results. 

First we have:

\begin{thm}\label{conv}
Let $V$ be a vertex operator algebra of central charge $c$
such that the associativity for intertwining operators
for $V$ (see \cite{H3}, \cite{H6} and \cite{H9}) hold, $W_{i}$
$V$-modules satisfying the 
$C_{2}$-cofiniteness condition,
$\tilde{W}_{i}$, $i=1, \dots, n$, $V$-modules, and
$\mathcal{Y}_{i}$, $i=1, \dots, n$, intertwining operators of 
types ${\tilde{W}_{i-1}\choose W_{i}\tilde{W}_{i}}$, respectively,
where, for convenience, we use the convention $\tilde{W}_{0}=
\tilde{W}_{n}$. 
Then for any $w_{i}\in W_{i}$, $i=1, \dots, n$,
the series
$$\tr_{\tilde{W}_{n}}\mathcal{Y}_{1}(w_{1}, z_{1})\cdots
\mathcal{Y}_{n}(w_{n}, z_{n})q^{L(0)-\frac{c}{24}}$$
is absolutely convergent in the region $1>|z_{1}|>\cdots 
>|z_{n}|>|q|>0$ and can be analytically extended
to a (multivalued) analytic function in the region 
$1>|z_{1}|, \dots, |z_{n}|>|q|>0$, 
$z_{i}\ne z_{j}$ when $i\ne j$. 
\end{thm}

Let $A_{j}$, $j\in \mathbb{Z}_{+}$, be complex numbers defined 
by 
\begin{eqnarray*}
\frac{1}{2\pi i}\log(1+2\pi i w)=\left(\exp\left(\sum_{j\in \mathbb{Z}_{+}}
A_{j}w^{j+1}\frac{\partial}{\partial w}\right)\right)w
\end{eqnarray*}
(see \cite{H5}).
Then we have
\begin{eqnarray*}
\lefteqn{\frac{1}{2\pi i}\log(1+e^{-2\pi iz}(w-z))}\nn
&&=\left(\frac{e^{-2\pi iz}}{2\pi i}\right)^{(w-z)
\frac{\partial}{\partial (w-z)}}
\left(\exp\left(\sum_{j\in \mathbb{Z}_{+}}
A_{j}(w-z)^{j+1}\frac{\partial}{\partial (w-z)}\right)\right)(w-z).
\end{eqnarray*}
Note that the composition inverse of $\frac{1}{2\pi i}\log(1+2\pi i w)$
is $\frac{1}{2\pi i}(e^{2\pi iw}-1)$ and thus we have
\begin{eqnarray*}
\frac{1}{2\pi i}(e^{2\pi iw}-1)=\left(\exp\left(-\sum_{j\in \mathbb{Z}_{+}}
A_{j}w^{j+1}\frac{\partial}{\partial w}\right)\right)w.
\end{eqnarray*}
Let $V$ be a vertex operator algebra. For any $V$-module 
$W$, we shall denote the operator 
$\sum_{j\in \mathbb{Z}_{+}}
A_{j}L(j)$
on $W$ by $L^{+}(A)$. Then 
$$e^{-\sum_{j\in \mathbb{Z}_{+}}
A_{j}L(j)}=e^{-L^{+}(A)}.$$

Let 
$W_{i}$, $\tilde{W}_{i}$  and $w_{i}\in W_{i}$ for 
$i=1, \dots, n$ be as in Theorem \ref{conv}.
For  any intertwining operators
$\mathcal{Y}_{i}$, $i=1, \dots, n$,  of 
types ${\tilde{W}_{i-1}\choose W_{i}\tilde{W}_{i}}$, respectively,
let
$$F_{\mathcal{Y}_{1}, \dots, \mathcal{Y}_{n}}(w_{1}, \dots, w_{n};
z_{1}, \dots, z_{n}; \tau)$$
be the analytic extension of 
\begin{eqnarray*}
\lefteqn{\tr_{\tilde{W}_{n}}\mathcal{Y}_{1}((2\pi ie^{2\pi iz_{1}})^{L(0)}
e^{-L^{+}(A)}
w_{1}, e^{2\pi iz_{1}})\cdots}\nn
&&\quad\quad\quad\quad\quad\quad\quad\quad\cdots
\mathcal{Y}_{n}((2\pi ie^{2\pi iz_{n}})^{L(0)}e^{-L^{+}(A)}
w_{n}, e^{2\pi iz_{n}})q^{L(0)-\frac{c}{24}}
\end{eqnarray*}
in the region $1>|e^{2\pi iz_{1}}|, \dots, |e^{2\pi iz_{n}}|>|q|>0$, 
$z_{i}\ne z_{j}$ when $i\ne j$, where $q=e^{2\pi i\tau}$.
We now consider the vector space $\mathcal{F}_{w_{1}, \dots, w_{n}}$
spanned by all such functions. Then we have the following result
obtained in \cite{H12}:

\begin{thm}
Let $V$ be a vertex operator algebra satisfying the following conditions:

\begin{enumerate}

\item For $n<0$, $V_{(n)}=0$  and $V_{(0)}=\mathbb{C}\mathbf{1}$.

\item Every $\mathbb{N}$-gradable $V$-module is completely reducible.

\item $V$ satisfies the $C_{2}$-cofiniteness condition.

\end{enumerate}

Then for any $V$-modules $W_{i}$ and $\tilde{W}_{i}$
and any
intertwining operators $\mathcal{Y}_{i}$   of 
types ${\tilde{W}_{i-1}\choose W_{i}\tilde{W}_{i}}$ ($i=1, \dots, n$), 
respectively,
and any 
$$\left(\begin{array}{cc}
\alpha&\beta\\
\gamma&\delta
\end{array}\right)\in SL(2, \mathbb{Z}),$$
we have
\begin{eqnarray*}
&{\displaystyle 
F_{\mathcal{Y}_{1}, \dots, \mathcal{Y}_{n}}
\left(\left(\frac{1}{\gamma\tau+\delta}\right)^{L(0)}w_{1}, \dots,
\left(\frac{1}{\gamma\tau+\delta}\right)^{L(0)}w_{n};
\frac{z_{1}}{\gamma\tau+\delta}, \dots, \frac{z_{n}}{\gamma\tau+\delta}; 
\frac{\alpha\tau+\beta}{\gamma\tau+\delta}\right)}&\nn
&\in \mathcal{F}_{w_{1}, \dots, w_{n}}&
\end{eqnarray*}

\end{thm}

We also need to show that the
traces corresponding to more general self-sewing of spheres 
also give modular 
invariants. This is related to the following problem 
on the uniformization of annuli:

\begin{prob}
Given an annulus, we know that it is conformally equivalent 
to a standard one.
Is it possible to obtain the conformal map from the 
standard one to the given one by sewing disks with analytically 
parametrized boundaries to pants successively? 
\end{prob}

In the problem, by obtaining a conformal
map by
sewing disks with analytically 
parametrized boundaries to pants we mean the following: 
Cut a hole in an annulus with analytically parametrized boundary
components, give the new boundary component 
an analytic parametrization  and  then sew a disk with 
analytically parametrized boundary to the resulting 
pants. The result is another annulus. If we sew 
such disks many times and choose the boundary parametrizations
such that the resulting annulus is conformally equivalent to the 
original one. Then we in fact obtain this conformal equivalence 
by sewing disks to pants successfully. 

We conjecture that this is possible, but probably one might need to sew
infinitely many times. If we do need to sew infinitely many times, then
we will need to take limits and will have to work with an analytic theory. 

Using all the results above in genus-one, we can construct 
genus-one weakly conformal field theories.

The third step is to construct weakly conformal field theories in any
genus. In particular, we need to construct locally convex completions of
irreducible modules for the vertex operator algebra. In fact, if
higher-genus correlation functions are constructed and their sewing
property is proved, one can construct the locally convex completions
from these correlation functions using the same method as the one used
in the construction of genus-zero completions in \cite{H8} and \cite{H10}
(see also \cite{CS} for some properties of the topology constructed
in \cite{H8}):
We first construct completions involving only genus-zero surfaces and
then add elements associated to higher-genus surfaces. Here by elements
associated to higher-genus surfaces, we mean the following: Consider a
Riemann surface with one positively oriented 
boundary component. If there is a conformal
field theory, this surface gives a map from $\mathbb{C}$ to $H$.  The
image of $1\in \mathbb{C}$ is called the {\it 
element associated to the surface}.
For example, in the genus-one case, for any $V$-module $W$, 
the one point function
$\tr_{W}Y(e^{2\pi izL(0)} e^{-L^{+}(A)} \cdot, e^{2\pi
iz})q^{L(0)-\frac{c}{24}}$ should be viewed as an element of $V^{*}$ and
should belong to the final completion of the module $V^{*}$.  In this
step, the main nontrivial problem is actually the convergence of
series corresponding to the sewing of higher-genus Riemann surfaces with
parametrized boundaries. The higher-genus modular invariance can be
obtained from the genus-zero duality (associativity) and the genus-one
modular invariance.  The details involve deep results in the analytic
theory of Teichm\"{u}ller spaces and moduli spaces and their generalizations
by Radnell in \cite{R} and have been
intensively studied by Radnell and the author.

The final step is to put two copies of a weakly conformal field theories
together suitably to get 
a conformal field theory. In \cite{S2}, Segal gave a description of 
how this step can be done for a 
weakly conformal field theories satisfying a certain unitarity condition.

Even though the problem of constructing conformal field theories has
not been worked out completely, the existing results already 
give many useful and
interesting results, for example, the Verlinde formula and the modular
tensor category structure mentioned above.  Here we discuss briefly the
problem of whether a vertex operator algebra or its suitable completion
is a double loop space.  We need the notion of operad which was
first formulated by May in \cite{May}. See for example, \cite{May},
\cite{KM} and \cite{H5} for the concepts needed in the theorem below,
including algebras over operads and the little disk operad.
The following theorem is obtained in \cite{H8} and \cite{H10} 
based on the geometry of vertex operator algebras in \cite{H1}, \cite{H2},
\cite{H5} and its reformulation in \cite{HL1} and \cite{HL2}
using the language of (partial) operads:

\begin{thm}
Let $V$ be a vertex operator algebra. 
Then there exists a unique minimal locally convex completion 
$H^{V}$ of $V$ such that $H^{V}$ is a genus-zero holomorphic 
conformal field theory. Here by minimal we mean that 
any genus-zero holomorphic conformal field theory containing 
$V$ must contain $H^{V}$. In particular, $H^{V}$ is an 
algebra over the little disk operad. 
\end{thm}

In \cite{May}, May proved the following recognition principle for double
loop spaces:

\begin{thm}
A space over the little disk operad has the weak homotopy type 
of a double loop space.
\end{thm}

Combining these two theorems, we conclude immediately that the locally
convex completion of a vertex operator algebra and the subset of its
nonzero elements have the weak homotopies of double loop spaces. But
actually the operad underlying the completion of a vertex operator
algebra is much richer than the little disk operad. So we can ask:

\begin{prob}
Can the operad underlying the completion of a vertex operator
algebra recognize more ``space-time'' properties and structures
of the completion of the algebra? 
Here by ``space-time'' properties and structures we
mean, for example, whether the algebra has a structure of 
a double loop space and so on. Especially, can it recognize
``space-time'' properties and structures homeomorphically, not only
(weak) homotopically? Can it recognize geometric properties and structures,
not only topological ones?
\end{prob}

An answer to the questions in this problem will undoubtedly provide
a deep understanding of the topological and geometric properties of 
conformal field theories.

\renewcommand{\theequation}{\thesection.\arabic{equation}}
\renewcommand{\thethm}{\thesection.\arabic{thm}}
\setcounter{equation}{0}
\setcounter{thm}{0}

\section{Open-closed conformal field theories}

Boundary conformal field theories \cite{C1} \cite{C2} are natural
generalizations of conformal field theories. In physics, boundary
conformal field theories describe realistic physical phenomena and
also describe $D$-branes in string theory. In this section, we
generalize Kontsevich's and Segal's notion of conformal field theory to
a notion of open-closed conformal field theory which incorporates the
axiomatic properties of  conformal field theories and boundary
conformal field theories.  In the topological case, a notion of
open-closed topological field theory has been introduced and studied
by Lazaroiu \cite{L} and by Moore and Segal \cite{Mo} \cite{S4}. 
Some part of the analogue in the conformal case
of the work \cite{L}, \cite{Mo} and \cite{S4}
has been done mathematically by Felder, Fr\"{o}hlich, 
Fuchs and Schweigert \cite{FFFS} and by Fuchs, 
Runkel and Schweigert \cite{FRS1} (see also \cite{FRS2}).
But in this case, the full analogues of (noncommutative) associative
algebras in the topological case discussed in \cite{L}, \cite{Mo}
and \cite{S4} should be algebras over the so-called "Swiss-cheese" operad.
These algebras still need to be constructed and studied
mathematically.
The goal of our
project on open-closed (or boundary)
conformal field theories is to establish the
corresponding formulations and results in the conformal case.
Here we shall discuss only the case that the surfaces involved 
are orientable. 
The unorientable theories can be formulated similarly using unorientable 
surfaces with conformal structures and the corresponding results 
can be obatined using double covering surfaces of the unoriented 
surfaces and certain operations in the state spaces. 
To avoid confusion, in the present section, we shall 
call a conformal field theory discussed in the preceding two sections
explicitly a closed conformal field theory.

Closed conformal field theories describe
closed string theory perturbatively.  So geometrically they are defined
in terms of Riemann surfaces with boundaries, where the boundary
components correspond to closed strings. In string theory, there are
also open strings. Also, $D$-branes are 
submanifolds of the space-time on which the end points of 
open strings move. To describe open strings
perturbatively, one still uses Riemann surfaces with boundaries, but
each boundary component is divided into two parts: the part
corresponding to open strings and the part corresponding to the motions
of the end points of open strings.

We first give a geometric  symmetric monoidal category 
$\mathcal{B}^{OC}$ for open-closed conformal field theories.
The object of $\mathcal{B}^{OC}$ are ordered sets 
of finitely many  copies of $[0, 1]$ and $S^{1}$.
The morphisms are conformal equivalence classes of Riemann surfaces with 
oriented, ordered and analytically parametrized segments of 
boundary components (parametrized by 
the copies of $[0, 1]$ in the objects) 
and oriented, ordered and analytically parametrized
boundary components (parametrized by the copies of $S^{1}$ in the 
objects). Note that by the definition of boundary component,
the boundary components 
containing segments parametrized
by $[0, 1]$ and the boundary compoents parametrized by $S^{1}$
are certainly disjoint. Also it is possible that there are
boundary components which are not parametrized at all, that is,
do not contain any segments parametrized by $[0, 1]$ and 
are also not parametrized by $S^{1}$. The monoidal structure
is also given by disjoint union as in the case of 
$\mathcal{B}$. Note that the category $\mathcal{B}$ discussed 
in Section 1 is a subcategory
of $\mathcal{B}^{OC}$. Also note  that for any functor or
projective functor $\Phi$ of monoidal categories from $\mathcal{B}^{OC}$ to
$\mathcal{T}$, if $\Phi(S^{1})=H^{C}$ and 
$\Phi([0, 1])=H^{O}$, then the images under $\Phi$ of 
objects of  must be in the tensor subcategory of $\mathcal{T}$
generated by $H^{C}$ and $H^{O}$.

An {\it open-closed conformal field theory} 
is a projective functor $\Phi$ from $\mathcal{B}^{OC}$ to the 
category $\mathcal{T}$  satisfying 
the following conditions: 
(i) Let $\Sigma$ be a morphism in $\mathcal{B}^{OC}$ from the 
ordered set of
$m$ copies of
$S^{1}$ and $p$  copies of $[0, 1]$ to the ordered set of 
$n$ copies of $S^{1}$
and $q$ copies of $[0, 1]$. Let $\Sigma_{\widehat{\widehat{i,j}}}$
(or $\Sigma_{\widehat{i,j}}$)
be the 
morphism from the ordered set of $m-1$  copies of $S^{1}$ and $p$ 
copies of $[0, 1]$ (or the ordered set of $m$ 
copies of $S^{1}$ and $p-1$ 
copies of $[0, 1]$)
to the ordered set of $n-1$ 
copies of $S^{1}$ and $q$  copies of $[0, 1]$
(or the ordered set of $n$ 
copies of $S^{1}$ and $q-1$  copies of $[0, 1]$)
obtained from $\Sigma$ by identifying the $i$-th copy of $S^{1}$ 
(or the $i$-th copy of $[0, 1]$)
in the domain of $\Sigma$ with the $j$-th copy of $S^{1}$ 
(or the $j$-th copy of $[0, 1]$)
in the codomain of $\Sigma$. Then the trace 
between the $i$-th copy of $H^{C}$ (or the $i$-th copy of $H^{O}$)
in the domain of $\Phi(\Sigma)$
and the $j$-th copy of $H^{C}$ (or the $j$-th copy of $H^{O}$)
 in the codomain of $\Phi(\Sigma)$ exists and is
equal to $\Phi(\Sigma_{\widehat{\widehat{i,j}}})$ (or 
$\Phi(\Sigma_{\widehat{i,j}})$). (ii)
Let $\Sigma$ be a morphism in $\mathcal{B}^{OC}$ from the ordered set of
$m$ copies of
$S^{1}$ and $p$ copies of $[0, 1]$ to the ordered set of 
$n$  copies of $S^{1}$
and $q$  copies of $[0, 1]$. Let $\Sigma_{i\to n+1}$
(or $\Sigma^{i\to n+1}$)
be the morphism from the ordered set of $m-1$  copies of
$S^{1}$ and $p$  copies of $[0, 1]$
(or the ordered set of $m$  copies of
$S^{1}$ and $p-1$  copies of $[0, 1]$) to the ordered set of 
$n+1$  copies of $S^{1}$ and $q$  copies of $[0, 1]$
(or the ordered set of 
$n$ copies of $S^{1}$ and $q+1$ copies of $[0, 1]$)
obtained from $\Sigma$ by changing the $i$-th
copy of $S^{1}$ (or the $i$-th copy of $[0, 1]$)
in the domain of $\Sigma$ 
to the $n+1$-st copy of $S^{1}$ (or the $q+1$-st copy 
of $[0, 1]$) in the codomain of $\Sigma_{i\to n+1}$
(or $\Sigma^{i\to n+1}$). Then $\Phi(\Sigma)$ and $\Phi(\Sigma_{i\to n+1})$
(or  $\Phi(\Sigma)$ and $\Phi(\Sigma_{i\to n+1})$)
are related by the map  obtained using the 
map $H^{C}\to (H^{C})^{*}$ ($H^{O}\to (H^{O})^{*}$) 
corresponding to the bilinear form 
$(\cdot, \cdot)_{C}$ (or $(\cdot, \cdot)_{O}$).

We consider an open-closed conformal field theory with the spaces
$H^{C}$ and $H^{O}$. Consider a Riemann surface $\Sigma_{C}^{O}$ with
two boundary components, one of its boundary component being negatively
oriented and parametrized by $S^{1}$ and the other boundary component
having one and only one positively oriented segment parametrized by $[0,
1]$.  Then by definition, the map corresponding to the equivalence class
of $\Sigma_{O}^{C}$ is a map $\iota_{\Sigma_{O}^{C}}$ from $H^{C}$ to $H^{O}$.
Similarly, associated to such a surface but with the opposite
orientations on the boundary components, denoted by $\Sigma_{O}^{C}$, we
have a map $\iota_{\Sigma_{O}^{C}}$ from $H^{O}$ to $H^{C}$. Note that
in the case of open-closed topological field theories studied in
\cite{L} and \cite{Mo} \cite{S4}, such maps are unique since all these surfaces are
topologically equivalent. But in an open-closed conformal field theory,
these maps depend on the equivalence classes of the surface
$\Sigma_{C}^{O}$ or $\Sigma_{O}^{C}$.  

It is also not difficult to see that the important Cardy condition
\cite{C2} is a consequence of the definition above. In fact, a cylinder
whose boundary are not parametrized at all 
can be obtained from a rectangle, with two opposite sides
parametrized by $[0, 1]$ and the other two sides not parametrized,
by identifying the two parametrized opposite sides. But
it can also be obtained by sewing three cylinders: The first has one 
boundary 
component which contains no
segments parametrized by $[0, 1]$ (so this component is not parametrized 
at all) and 
has one positively
oriented boundary component parametrized by $S^{1}$; the second one has
one negatively oriented and one positively oriented boundary components
which are both parametrized by $S^{1}$; the third has one boundary component 
which contains no segments parametrized by $[0, 1]$ (so as in the case
of the first cylinder, this componenet is not parametrized at all)
and has one negatively oriented
boundary component parametrized by $S^{1}$. According to the axioms,
the map corresponding to the cylinder obtained from a rectangle by 
identifying two opposite sides is the trace of the map 
corresponding to the rectangle. According to the axioms again,
the other way of obtaining this cylinder means that it is 
in fact the composition of three maps: The first is a map from 
$\mathbb{C}$ to $H^{C}$; the second is a map from $H^{C}$ to $H^{C}$;
the third is a map from $H^{C}$ to $\mathbb{C}$. Since the first and
the last map are actually equivalent to an element of $H^{C}$
and an element of $(H^{C})^{*}$, respectively, we see that 
the composition of the three maps above is equivalent to the 
matrix element 
between these two elements of $H^{C}$ and $(H^{C})^{*}$
of the map from $H^{C}$ to $H^{C}$. Now by the axioms again, we see
that the trace of the map 
corresponding to the rectangle is equal to this matrix element. 
This is exactly the {\it Cardy condition}.

As in the case of closed conformal field theories,
the most important problem is the following:

\begin{prob}
Construct open-closed conformal field theories satisfying the definition 
above.
\end{prob}

From the definition, we see that open-closed conformal field theories
must have closed conformal field theories as subtheories. Therefore it
is natural to try to construct open-closed conformal field theories from
closed conformal field theories. Here is a strategy used by physicists: 
A Riemann surface for
open-closed conformal field theories can be doubled to obtain a Riemann
surface for closed conformal field theories. This procedure establishes a
connection, at least geometrically, between closed 
conformal field theories and
open-closed conformal field theories. This connection provides a
concrete way to construct open-closed conformal field theories from
closed conformal field theories.

For closed conformal field theories, one starts with the genus-zero case.
For open-closed conformal field theories, the construction 
of genus-zero theories 
is again the first step.
The simplest genus-zero surfaces in this case are
Riemann surfaces with only one 
boundary component and with segments of the boundary component 
parametrized by $[0, 1]$.
In particular, these surfaces do not have 
boundary components parametrized by $S^{1}$.

We consider certain such special Riemann surfaces.  They are the closed
upper half unit disk with some smaller disks inside and 
some smaller upper half disks centered on the
real line removed and with the obvious parametrizations for all the full 
and half
circles.  These Riemann surfaces form an operad and is called the
``Swiss cheese'' operad by Voronov \cite{Vo}. Just as in the case of
closed conformal field theories, our first step is to construct algebras
over this operad. In fact, in an open-closed topological field theory
discussed in \cite{L} and \cite{Mo}, there is an associative algebra
which in general is not commutative.  Algebras over the Swiss cheese
operad should be viewed as analogues of such algebras in open-closed
conformal field theories.  

To construct algebras over the Swiss cheese operad, we form doubles of
these surfaces.  These doubles are disks with smaller disks removed and
are special elements of the little disk operad. So any algebra over the
little disk operad is automatically an algebra over the Swiss cheese
operad.  In particular, the locally convex completion of a vertex
operator algebra is an algebra over the Swiss cheese operad. But there
are more such algebras.  First, the vertex operator algebras
corresponding to the upper halves and to the lower halves of the doubles
can differ by an isomorphism. Second, since the middle disks are always
centered on the real line, we do not have to worry about the multivalued
property of the corresponding operators and thus we can place elements
of modules for the vertex operator algebras at these middle disks.  In
particular, we can construct algebras over the Swiss cheese operad from
suitable subalgebras of the intertwining operator algebras constructed
from modules and intertwining operators for vertex operator
algebras. This method of constructing algebras over the Swiss cheese
operad is studied by Kong and the author in \cite{HK}. 
In particular, in \cite{HK}, a notion of open-string vertex algebra is
introduced and studied and 
examples of such algebras are constructed. It is established in \cite{HK}
that the category of open-string vertex algebras
is equivalent to the category of so-called differentiable-meromorphic
algebras 
over a partial operad which is an extension of the swiss cheese operad.

Note that in this case the vertex operator algebras corresponding to the
upper and lower half disks can differ by an isomorphism.  Thus we have
to make sure that they agree on the real line and this requirement is
called a {\it boundary condition}.  In general it is certainly not true
that these boundary conditions for vertex operators hold when the vertex
operators act on all elements of a module for the algebra.  But it is
possible that when the vertex operators act on certain special elements
of modules for the algebra, a boundary condition is satisfied.  These
special elements are called {\it boundary states} for the given boundary
condition. In all known examples in physics (certainly still to be
constructed mathematically), 
the Hilbert space $H^{O}$ for the
open-closed conformal field theory can be decomposed as a direct sum of
spaces $H^{O}_{ab}$ which are suitable completions of the spaces of
boundary states obtained by imposing boundary conditions labeled by $a$
and $b$ at the boundary segments corresponding to the movements of the
two end points $0$ and $1$, respectively, of the open string described
by $[0, 1]$.  Note that for boundary states, it is important to have
locally convex completion of modules for the vertex operator algebra
since boundary states are in general infinite series of elements of the
homogeneous subspaces of modules. A boundary condition at the boundary
segment corresponding to the movement of an end point of an open string
amounts to exactly a so-called {\it $D$-brane} in the string theory.

In analogy with the problem of whether a vertex operator algebra or its
completion is a double loop space, for open-closed conformal field
theories, we also have the analogous problem of whether an algebra or a
space over the Swiss cheese operad has a certain loop space structure.
Given a topological space with a base point, consider paths from the base
point to arbitrary points. Now consider based
loops in the space of all such
paths. It is easy to see that the Swiss cheese operad acts on the space
of such loops. 

\begin{prob}
Can we recognize such a loop space structure from a structure 
of a space over the Swiss cheese operad?
\end{prob}

If possible, then we have a topological space such that open strings 
(paths) move in it. Such a picture might be helpful 
in the study of $D$-branes in a genegral conformal field theory
background since they can then be viewed as 
subsets in the topological space.

\noindent {\small \sc Department of Mathematics, Rutgers University,
110 Frelinghuysen Rd., Piscataway, NJ 08854-8019}

\vskip 1em

\noindent {\em E-mail address}: yzhuang@math.rutgers.edu

\end{document}